\documentclass[12pt,english]{amsart}
\usepackage[T1]{fontenc}
\usepackage[latin9]{inputenc}
\usepackage{geometry}
\usepackage{amstext}
\usepackage{amsthm}
\usepackage{amssymb}
\usepackage{ulem}
\usepackage{color}

\makeatletter
\numberwithin{equation}{section}
\numberwithin{figure}{section}
\theoremstyle{plain}
\newtheorem{thm}{\protect\theoremname}
\theoremstyle{plain}

\usepackage{fancyhdr}
\usepackage{cite}

\newcommand{\cC}{\mathcal{C}}

\newcommand{\cM}{\mathcal{M}}

\newcommand{\cO}{\mathcal{O}}

\newcommand{\abs}[1]{\ensuremath{|#1|}}

\newcommand{\Abs}[1]{\ensuremath{\left|#1\right|}}

\newcommand{\Norm}[2]{\ensuremath{\left|\!\left|#1\right|\!\right|_{#2}}}

\makeatother

\usepackage{babel}
\providecommand{\propositionname}{Proposition}
\providecommand{\theoremname}{Theorem}

\begin{document}

\title{Sharp estimate on the resolvent of a finite-dimensional contraction}

\author{Karine Fouchet }

\address{Karine Fouchet, Institut de Mathematiques, UMR 7373, Aix-Marseille
Universite, 39 rue F. Joliot Curie, 13453 Marseille Cedex 13, France}

\email{karine.isambard@univ-amu.fr}



\thanks{The work is supported by the project ANR 18-CE40-0035}
\begin{abstract}
We compute an asymptotic formula for the supremum of the resolvent
norm $\Norm{(\zeta-T)^{-1}}{}$ over $\left|\zeta\right|\geq1$ and
contractions $T$ acting on an $n-$dimensional Hilbert space, whose
spectral radius does not exceed a given $r\in(0,\,1)$. We prove that
this supremum is achieved on the unit circle by an analytic Toeplitz
matrix. 
\end{abstract}

\maketitle

\section{\label{sec:Statement-of-the}Introduction}

\subsection{Notation and main result }

Let $\cM_{n}$ be the set of complex $n\times n$ matrices and let
$\Norm{T}{}$ denote the operator norm of $T\in\cM_{n}$ associated
with the Hilbert norm on $\mathbb{C}^{n}$. We denote by $\sigma=\sigma(T)$
the spectrum of $T$ and by $\rho(T)=\max_{\lambda\in\sigma(T)}\left|\lambda\right|$
its spectral radius. In our discussion we will assume that $\Norm{T}{}\leq1$
and call such $T$ a \textit{contraction}. Let $\cC_{n}\subset\cM_{n}$
denote the set of all contractions. We denote by $R(\zeta,\,T)=(\zeta-T)^{-1}$
the resolvent of $T$ at point $\zeta\,\notin\sigma$ and are interested
in sharp estimates on $\Norm{R(\zeta,\,T)}{}$ in terms of the spectral
data of $T$. E. B. Davies and B. Simon~\cite{DaSi} proved that
the supremum
\[
\sup_{\abs{\zeta}\geq1}\sup\left\{ {d}(\zeta,\,\sigma)\Norm{R(\zeta,\,T)}{}\ :\ T\in\cC_{n}\right\} 
\]
 is attained for $\abs{\zeta}=1$ and that for such $\zeta$ we have
\begin{equation}
\sup\left\{ {d}(\zeta,\,\sigma)\Norm{R(\zeta,\,T)}{}\ :\ T\in\cC_{n}\right\} =\cot(\frac{\pi}{4n}).\label{eq:DS_result}
\end{equation}
In this paper we will show \textendash{} using different methods from
those in ~\cite{DaSi} \textendash{} that given $r\in(0,1)$ the
supremum
\[
\sup_{\left|\zeta\right|\geq1}\mbox{sup}\left\{ \Norm{R(\zeta,\,T)}{}\ :\ T\in\cC_{n},\,\rho(T)\leq r\right\} 
\]
is again attained for $\abs{\zeta}=1$ and that for such $\zeta$
we have
\begin{equation}
\mathcal{R}_{n,\,r}:=\mbox{sup}\left\{ \Norm{R(\zeta,\,T)}{}\ :\ T\in\cC_{n},\,\rho(T)\leq r\right\} \sim\frac{2}{\pi}\frac{1+r}{1-r}n\label{eq:Our_result}
\end{equation}
as $n\rightarrow\infty.$ The paper is organized as follows: Section
\ref{subsec:Motivation} motivates our study relating our question
to problems previously considered by V. Pt\'ak and N. Young \cite{Pta,PtYo}.
In Section \ref{subsec:Statement-of-our} we state our results
in full detail (see Theorems \ref{thm:Ptak_type_result} - 2 -  \ref{prop} below) and exhibit an analytic $n\times n$ Toeplitz matrix
that achieves the supremum in (\ref{eq:Our_result}). We also mention
the techniques from matrix analysis on the one hand and from interpolation
theory on the other hand, which we employ to prove Theorem \ref{thm:Ptak_type_result}, Theorem 2
and Theorem \ref{prop} respectively. In Section \ref{sec:Main-ingredients}
we first lay down the required definitions and basic results from
the theory of model spaces and their operators, to the approach we
choose to estimate $\Norm{R(\zeta,\,T)}{}$. Then we relate our question
to interpolation theory by expressing the supremum in (\ref{eq:Our_result})
as an interpolation quantity in the algebra $H^{\infty}$ of bounded
holomorphic functions in the unit disc $\mathbb{D}=\{z\in\mathbb{C}:\,\vert z\vert<1\}$.
Finally Section \ref{sec:Proofs} collects the proofs of Theorems~\ref{thm:Ptak_type_result}--\ref{prop}. 

\subsection{\label{subsec:Motivation}Motivation}

Our motivation to consider such of a variant of Davies-Simon's question
comes from problems studied by Pt\'ak and Young \cite{Pta,PtYo}. For
a finite sequence $\sigma$ in $\mathbb{D}$, we denote by $P_{\sigma}$
the monic polynomial with zero set $\sigma$ (counted with multiplicities).
Given a finite sequence $\sigma$ in $\mathbb{D}$ and $f\in\mathcal{H}ol(\mathbb{D}),$
$\mathcal{H}ol\left(\mathbb{D}\right)$ being the space of holomorphic
functions in $\mathbb{D}$, we introduce the quantity 
\[
\mathcal{S}(f,\,\sigma)=\mbox{sup}\left\{ \Norm{f(T)}{}:\,T\in\cC_{n},\,m_{T}=P_{\sigma}\right\} ,
\]
where $m_{T}$ stands for the minimal polynomial of $T\in\cC_{n}$.
The case $f\vert_{\sigma}={\displaystyle z^{k}\vert_{\sigma}}$ (estimates
on the norm of the powers of an $n\times n$ contraction) is considered
by Pt\'ak in \cite{Pta} who proved that the quantity 
\[
\mbox{sup}\left\{ \Norm{T^{k}}{}:\:T\in\cC_{n},\,\rho(T)\leq r\right\} ,
\]
is achieved by an $n\times n$ analytic Toeplitz matrix. In the same
time Young \cite{You} proved among other things that at fixed $r$
the above quantity is less than 
\[
nr+O(r^{2}).
\]
Later Pt\'ak and Young \cite[Section 2]{PtYo} considered the more general
quantity 
\[
\sup\big\{\mathcal{S}(f,\,\sigma):\:\max_{\lambda\in\sigma}\left|\lambda\right|\leq r\},
\]
where $\mathcal{S}(f,\,\sigma)$ is defined above, $\sigma$ is a
finite sequence of $\mathbb{D}$ and $f$ is an analytic polynomial.
They proved \cite[Section 2, p. 365]{PtYo} that the above supremum
is achieved by an $n\times n$ analytic Toeplitz matrix. Here we study
the case $f\vert_{\sigma}=(\zeta-z)^{-1}\vert_{\sigma}$ where $\abs{\zeta}=1$
(estimates on the norm of the resolvent of an $n\times n$ contraction)
and our estimate (\ref{eq:Our_result}) can be reformulated using
Pt\'ak and Young's notation as 
\[
\sup\big\{\mathcal{S}\left(\left(\zeta-z\right)^{-1},\,\sigma\right):\:\max_{\lambda\in\sigma}\left|\lambda\right|\leq r\}\sim\frac{2}{\pi}\frac{1+r}{1-r}n,\qquad n\rightarrow\infty.
\]

\subsection{\label{subsec:Statement-of-our}Statement of our result and main
ingredients to its proof}

We prove the following theorem for which we recall that the definition
of $\mathcal{R}_{n,\,r}$ is given by the left-hand side of (\ref{eq:Our_result}). 

\begin{thm}
\label{thm:Ptak_type_result} Given $n\geq1$, $r\in(0,\,1)$ $\zeta\in\mathbb{C}\setminus\mathbb{D}$
and $T\in\cC_{n}$ such that $\rho(T)\leq r$ we have
\begin{equation}
\Norm{R(\zeta,\,T)}{}\leq \Norm{\big(1-T^{*}\big)^{-1}}{}=\frac{1}{1-r}\Norm{X_{1+r}}{},\label{eq:H_3}
\end{equation}
where $T^{*} \in \cC_{n} $ is the analytic $n\times n$ Toeplitz matrix 
\[
T^{*}=\left(\begin{array}{ccccc}
r & 0 & \ldots & \ldots & 0\\
1-r^{2} & r & \ddots & . & \vdots\\
-r(1-r^{2}) & 1-r^{2} & r & \ddots & \vdots\\
\vdots & \ddots & \ddots & \ddots & 0\\
(-r)^{n-2}(1-r^{2}) & \ldots & -r(1-r^{2}) & 1-r^{2} & r
\end{array}\right)
\]
and the analytic $n\times n$ Toeplitz matrix $X_{1+r}$ is entry-wise
given by:
\[
\left(X_{1+r}\right)_{ij}=\Bigg\{\begin{array}{ll}
0 & \mbox{if }i<j\\
1 & \mbox{if }i=j\\
1+r & \mbox{if }i>j
\end{array}.
\]
In particular, we have $\mathcal{R}_{n,\,r} = \Norm{\big(1-T^{*}\big)^{-1}}{}$.
\end{thm}
\begin{thm}
At fixed $r\in(0,\,1)$ and in the limit of large $n$ we
have
\begin{equation}
\mathcal{R}_{n,\,r}\sim\frac{2}{\pi}\frac{1+r}{1-r}n.\label{eq:H_2}
\end{equation}
\end{thm}

The proof of (\ref{eq:H_3}) is based on an application of the Commutant
Lifting Theorem of B. Sz.-Nagy and C. Foia\c{s} \cite{NaFo,FoFr,Sar}
to the rational function $f(z)=\frac{1}{\zeta-z}$. Our approach to
compute the norm of $\Norm{X_{1+r}}{}$ and prove the asymptotic formula
(\ref{eq:H_2}), follows the one from \cite{Sze,SzZa} and builds
on the techniques developed in \cite{Ege}.\\

Using a purely interpolation-theoretic approach we obtain the following
results, which are asymptotically weaker than formula (\ref{eq:H_2}),
but hold for any $n\geq1$ and $r\in[0,1)$.

\begin{thm}
\label{prop}Let $T\in\cC_{n}$ with minimal polynomial $m$ and spectrum
$\sigma=\left(\lambda_{1},\lambda_{2},\,\dots,\,\lambda_{\abs{m}}\right)\in\mathbb{D}^{\abs{m}},$
where $\abs{m}$ denotes the degree of $m$. We have
\begin{equation}
\Norm{R(\zeta,\,T)}{}\leq\sum_{k=1}^{\abs{m}}\frac{1+\vert\lambda_{k}\vert}{1-\vert\lambda_{k}\vert},\qquad\abs{\zeta}\geq1.\label{eq:ptak_upper_bd}
\end{equation}
In particular, the upper estimate 
\[
\mathcal{R}_{n,\,r}\leq\frac{1+r}{1-r}n,
\]
holds for any $n\geq1$ and $r\in[0,1)$. Moreover the lower estimate
\begin{equation}
\mathcal{R}_{n,\,r}\geq\frac{1}{2}\frac{1+r}{1-r}n+\frac{1}{2},\label{eq:weaker_ptak}
\end{equation}
holds also for any $n\geq1$ and $r\in[0,1)$.
\end{thm}

The asymptotically sharp prefactor $\frac{2}{\pi}$ does not appear
in the upper estimate on $\mathcal{R}_{n,\,r}$ from Theorem \ref{prop},
which is therefore asymptotically weaker than formula (\ref{eq:H_2}).
The proof of (\ref{eq:ptak_upper_bd}) is obtained by applying Von Neumann's
inequality to $T\in\cC_{n}$ and to the orthogonal projection of the
function $z\mapsto\frac{1}{\zeta-z}$ onto $K_{B_{\sigma}},$ $\sigma$
being the spectrum of $T.$ The lower estimate (\ref{eq:weaker_ptak})
is derived by making use of $H^{\infty}-$interpolation techniques developed
in \cite{BaZa1,Zar}.

\smallskip{}

\section{\label{sec:Main-ingredients}Main ingredients}

This section lays down the required definitions and basic results used in
the approach we choose to estimate $\Norm{R(\zeta,\,T)}{}$. We
begin with the definition of model spaces and model operators.

\subsection{\label{subsec:Model-spaces-and}Model spaces and model operators}

We denote by $H^{\infty}$ the Banach algebra of bounded analytic
functions in $\mathbb{D}$, endowed with the supremum norm on $\mathbb{D}$:
$\Norm{f}{\infty}=\sup_{z\in\mathbb{D}}\left|f(z)\right|$. We recall
that the standard Hardy space $H^{2}$ is defined as the subspace
of $\mathcal{H}ol(\mathbb{D})$ consisting of those functions $f$
such that 
\[
\Vert f\Vert_{H^{2}}^{2}:=\sup_{0\leq r<1}\int_{\mathbb{\partial\mathbb{D}}}\left|f(rz)\right|^{2}{\rm d}m(z)<\infty,
\]
where $m$ is the normalized Lebesgue measure on the unit circle $\partial\mathbb{D}:=\left\{ z\in\mathbb{C}:\,|z|=1\right\} $.
Endowed with $\Vert\cdot\Vert_{H^{2}}$, $H^{2}$ is a Hilbert space.

Let $\sigma$ be a finite sequence of points in $\mathbb{D}.$ The
finite Blaschke product $B=B_{\sigma}$ corresponding to $\sigma$
is defined by 
\[
B=B_{\sigma}=\prod_{\lambda\in\sigma}b_{\lambda},
\]
where $b_{\lambda}(z)=\frac{\lambda-z}{1-\overline{\lambda}z}$ is the
Blaschke factor corresponding to $\lambda\in\mathbb{D}$. Then one
defines the model space $K_{B}$ as the finite dimensional subspace
of $H^{2}$ given by 
\[
K_{B}=\left(BH^{2}\right)^{\perp}=H^{2}\ominus BH^{2}.
\]
The reproducing kernel of the model space $K_{B}$ corresponding to
a point $\zeta\in\mathbb{D}$ is of the form 
\[
k_{\zeta}^{B}(z)=\frac{1-\overline{B(\zeta)}B(z)}{1-\overline{\zeta}z}=(1-\overline{B(\zeta)}B(z))k_{\zeta}(z),
\]
where $k_{\zeta}(z)=\frac{1}{1-\bar{\zeta}z}$ is the Cauchy kernel at
$\zeta.$ If $\sigma=(\lambda_{1},\dots,\lambda_{n})\in\mathbb{D}^{n}$,
then we set  $f_{j}(z)=k_{\lambda_{j}}(z)=\frac{1}{1-\overline{\lambda_{j}}z}$,
$j=1,\ldots,n$. Observe that $\Vert f_{j}\Vert_{H^{2}}=\left(1-\vert\lambda_{j}\vert^{2}\right)^{-1/2}$.
Now one can check that the family $\left(e_{k}\right)_{1\leq k\leq n}$
given by 
\[
e_{1}=\frac{f_{1}}{\Vert f_{1}\Vert_{H^{2}}}\,\quad\mbox{and}\quad e_{k}={\displaystyle \prod_{j=1}^{k-1}}b_{\lambda_{j}}\frac{f_{k}}{\Vert f_{k}\Vert_{H^{2}}},\quad k=2,\ldots,n,
\]
is an orthonormal basis of $K_{B}$ (known as the Takenaka-Malmquist-Walsh
basis, see \cite[p. 117]{Nik3}). In particular for any $\zeta\in\overline{\mathbb{D}}$
we have:
\begin{equation}
k_{\zeta}^{B}=P_{B}\left(k_{\zeta}\right)=\sum_{k=1}^{n}\overline{e_{k}(\zeta)}e_{k}.\label{eq:Taylor_Malmq_k_zeta}
\end{equation}

Further, we define the model operator $M_{B}$ acting on $K_{B}$
as follows:

\[
M_{B}:\left\{ \begin{array}{ccc}
K_{B} & \rightarrow & K_{B}\\
f & \mapsto & P_{B}(zf),
\end{array}\right.
\]
where $P_{B}$ denotes the orthogonal projection on $K_{B}.$ We finally
denote by $\hat{M_{\sigma}}$ the matrix representation of $M_{B}$
with respect to the Takenaka-Malmquist-Walsh basis $(e_{k})_{1\leq k\leq n}$
of $K_{B}.$ By \cite[Proposition III.4]{Sze},

\[
\left(\hat{M_{\sigma}}\right)_{ij}=\left\{ \begin{array}{ll}
0 & \text{if }i<j\\
\lambda_{i} & \text{if }i=j\\
(1-\abs{\lambda_{i}}^{2})^{1/2}(1-\abs{\lambda_{j}}^{2})^{1/2}\prod_{\mu=j+1}^{i-1}\left(-\bar{\lambda}_{\mu}\right) & \text{if }i>j,
\end{array}\right.
\]
where $\left(\hat{M_{\sigma}}\right)_{ij}$ stands for the $i,j$
entry of $\hat{M_{\sigma}}$ , and the empty product is defined to
be 1.\smallskip{}

We refer the reader to \cite{Nik2,Nik3} for a thorough study of model
operators and model spaces.

\subsection{\label{subsec:Link-to-interpolation}Link to interpolation theory}

Given a Blaschke sequence $\sigma$ in $\mathbb{D}$ and $f\in H^{\infty}$
it is possible to evaluate $\mathcal{S}(f,\,\sigma)$ as follows:
\begin{equation}
\mathcal{S}(f,\,\sigma)=\Norm{f}{H^{\infty}/B_{\sigma}H^{\infty}}=\Norm{f(M_{B_{\sigma}})}{},\label{eq:Sarason}
\end{equation}
where $M_{B_{\sigma}}$ is the compression of the multiplication operation
by $z$ to the model space $K_{B_{\sigma}}$. This formula is due
to N. K. Nikolski \cite[Theorem 3.4]{Nik1} while the last equality
is a well-known corollary of Commutant Lifting Theorem of B. Sz.-Nagy
and C. Foia\c{s} \cite{NaFo,FoFr,Sar}. \\
It is shown in \cite{BaZa1} that the above equality on $\mathcal{S}(f,\,\sigma)$
naturally extends to any $f\in\mathcal{H}ol(\mathbb{D})$ as follows.
We reproduce the proof of this fact for completeness. There exists
an analytic polynomial $p$ interpolating $f$ on the finite set $\sigma$.
Therefore for any $T\in\cC_{n}$ with $m_{T}=P_{\sigma}$ and $\sigma\subset\mathbb{D}$,
we have $f(T)=p(T)$ (since $f=p+m_{T}h$ for some $h\in\mathcal{H}ol(\mathbb{D})$).
Hence, 
\begin{eqnarray*}
\mathcal{S}(f,\,\sigma) & = & \mathcal{S}(p,\,\sigma)=\Norm{p}{H^{\infty}/B_{\sigma}H^{\infty}}\\
 & = & \Norm{p(M_{B_{\sigma}})}{}=\Norm{f(M_{B_{\sigma}})}{}.
\end{eqnarray*}
Here we used (\ref{eq:Sarason}) applied to $p$. Moreover 
\begin{align*}
\Norm{p}{H^{\infty}/B_{\sigma}H^{\infty}} & =\inf\{\Norm{p+B_{\sigma}h}{\infty}:\:h\in H^{\infty}\}\\
 & =\inf\{\Norm{g}{\infty}:\:\:g\vert_{\sigma}=p\vert_{\sigma},\:g\in H^{\infty}\}\\
 & =\inf\{\Norm{g}{\infty}:\:\:g\vert_{\sigma}=f\vert_{\sigma},\:g\in H^{\infty}\}.
\end{align*}
We conclude that 
\begin{equation}
\mathcal{S}(f,\,\sigma)=\inf\,\big\{\Norm{g}{\infty}:\:g\in H^{\infty},\:g\vert_{\sigma}=f\vert_{\sigma}\big\}.\label{eq:interp_identity}
\end{equation}
In particular we will apply (\ref{eq:interp_identity}) with $f(z)=(\zeta-z)^{-1}$
to prove the lower estimate (\ref{eq:weaker_ptak}) from Theorem
\ref{prop}.

\section{\label{sec:Proofs}Proofs}

\subsection{{Proof of the upper bound in Theorem \ref{prop}}}
\begin{proof}
Let $T\in\cC_{n}$ with minimal polynomial $m$ and spectrum $\sigma$.
Let $B=B_{\sigma}$ the finite Blaschke product corresponding to $\sigma=\left(\lambda_{1},\lambda_{2},\,\dots,\,\lambda_{\abs{m}}\right)$.
We suppose that the spectral radius of $T$ is less than $r\in(0,\,1)$.
We first assume that $\vert\zeta\vert>1$ (the more general case $\zeta\in\overline{\mathbb{D}}$
will follow immediately by continuously moving $\zeta$ towards the
boundary $\partial\mathbb{D}$ of $\mathbb{D}$) and consider the
orthogonal projection of $\frac{1}{\zeta-z}$ onto the model space
$K_{B}:$ 
\[
g(z)=\frac{1}{\zeta}P_{B}\big(k_{1/\bar{\zeta}})(z).
\]
The function $g$ clearly interpolates $\frac{1}{\zeta-z}$ on $\sigma.$
Moreover expanding $g$ on the Takena-Malmquist-Walsh $(e_{k})_{k=1}^{\abs{m}}$
basis of $K_{B}$ we get
\begin{align*}
g=\frac{1}{\zeta}P_{B}\big(k_{1/\bar{\zeta}}) & =\frac{1}{\zeta}\sum_{k=1}^{\abs{m}}\left\langle k_{1/\bar{\zeta}},\,e_{k}\right\rangle e_{k}\\
 & =\frac{1}{\zeta}\sum_{k=1}^{\abs{m}}\overline{e_{k}(1/\bar{\zeta})}e_{k}.
\end{align*}
 That is to say that for any $u\in\partial\mathbb{D}$ 
\[
g(u)=\sum_{k=1}^{\abs{m}}\frac{1-\vert\lambda_{k}\vert^{2}}{\zeta-\lambda_{k}}\overline{\prod_{j=1}^{k-1}b_{\lambda_{j}}(1/\bar{\zeta})}\prod_{j=1}^{k-1}b_{\lambda_{j}}(u)\frac{1}{1-\overline{\lambda_{k}}u}.
\]
Clearly we have
\[
\Norm{g}{\infty}\leq\sum_{k=1}^{\abs{m}}\frac{1+\vert\lambda_{k}\vert}{1-\vert\lambda_{k}\vert},
\]
and the proof of (\ref{eq:ptak_upper_bd}) follows by making use of
Von Neumann's inequality:
\[
\Norm{R(\zeta,\,T)}{}=\Norm{g(T)}{}\leq\Norm{g}{\infty}.
\]
\end{proof}

\subsection{\label{subsec:The-lower-bounds}{Proof of the lower bound in Theorem \ref{prop}}}

In this paragraph we prove (\ref{eq:weaker_ptak}) by using $H^{\infty}-$interpolation
techniques. 
\begin{proof}
Again we first assume that $\vert\zeta\vert>1$ (the more general
case $\zeta\in\overline{\mathbb{D}}$ will follow immediately by continuously
moving $\zeta$ towards $\partial\mathbb{D}$). Consider $\lambda\in\mathbb{D}$,
$B=b_{\lambda}^{n}$ and the function 
\begin{align*}
g(z) & =\frac{1}{\zeta}P_{B}\big(k_{1/\bar{\zeta}})(z)\\
 & =\frac{1}{\zeta}\frac{1-\overline{B(1/\bar{\zeta})}B(z)}{1-z/\zeta}\\
 & =\frac{1}{\zeta}\sum_{k=0}^{n-1}(\overline{b_{\lambda}(1/\bar{\zeta})})^{k}\frac{1}{1-\frac{\lambda}{\zeta}}(b_{\lambda}(z))^{k}\frac{1-\vert\lambda\vert^{2}}{1-\bar{\lambda}z},
\end{align*}
where the last equality is due to (\ref{eq:Taylor_Malmq_k_zeta}).
Clearly, 
\[
g-\frac{1}{\zeta} k_{1/\bar{\zeta}}\in B\mathcal{H}ol(\mathbb{D}),
\]
and in particular 
\[
\Norm{g}{H^{\infty}/b_{\lambda}^{n}H^{\infty}}=\Norm{g(M_{B})}{}=\Norm{(\zeta-M_{B})^{-1}}{}.
\]
Since $1-\overline{\lambda}b_{\lambda}(z)=\frac{1-\bar{\lambda}z-\vert\lambda\vert^{2}+\bar{\lambda}z}{1-\bar{\lambda}z}=\frac{1-\vert\lambda\vert^{2}}{1-\bar{\lambda}z}$ and $b_\lambda$ is self-inverse,
we have 
\begin{align*}
g\circ b_{\lambda}(z) & =\sum_{k=0}^{n-1}(\overline{b_{\lambda}(1/\bar{\zeta})})^{k}\frac{1}{\zeta-\lambda}z^{k}\frac{1-\vert\lambda\vert^{2}}{1-\overline{\lambda}b_{\lambda}(z)}\\
 & =\frac{1}{\zeta-\lambda}\sum_{k=0}^{n-1}(\overline{b_{\lambda}(1/\bar{\zeta})})^{k}z^{k}\left(1-\bar{\lambda}z\right)\\
 & =\frac{1}{\zeta-\lambda}\left(1-\bar{\lambda}(\overline{b_{\lambda}(1/\bar{\zeta})})^{n-1}z^{n}+\sum_{k=1}^{n-1}(\overline{b_{\lambda}(1/\bar{\zeta})})^{k}z^{k}-\bar{\lambda}\sum_{k=1}^{n-1}(\overline{b_{\lambda}(1/\bar{\zeta})})^{k-1}z^{k}\right)\\
 & =\frac{1}{\zeta-\lambda}\left(1-\bar{\lambda}(\overline{b_{\lambda}(1/\bar{\zeta})})^{n-1}z^{n}+\left(\overline{b_{\lambda}(1/\bar{\zeta})}-\bar{\lambda}\right)\sum_{k=1}^{n-1}(\overline{b_{\lambda}(1/\bar{\zeta})})^{k-1}z^{k}\right)\\
 & =\frac{1}{\zeta-\lambda}\left(1-\bar{\lambda}(\overline{b_{\lambda}(1/\bar{\zeta})})^{n-1}z^{n}-\frac{1-\vert\lambda\vert^{2}}{\zeta-\lambda}\sum_{k=1}^{n-1}(\overline{b_{\lambda}(1/\bar{\zeta})})^{k-1}z^{k}\right)
\end{align*}
Given $r\in(0,\,1)$ and $-\frac{1}{r}<\zeta<-1,$ we consider $\lambda=-r.$
The function 
\[
\Psi_{n}(z):=\,g\circ b_{-r}(z)=\frac{1}{\zeta+r}\bigg(1+r\left(-\frac{1+r\zeta}{\zeta+r}\right)^{n-1}z^{n}-\frac{1-r^{2}}{r+\zeta}\sum_{k=1}^{n-1}\left(-\frac{1+r\zeta}{r+\zeta}\right)^{k-1}z^{k}\bigg),
\]
is an analytic polynomial of degree $n$.
We need a lower estimate for $\Vert g\Vert_{H^{\infty}/b_{\lambda}^{n}H^{\infty}}$:
to this end consider $G$ such that $R_{n}-G\in b_{\lambda}^{n}\mathcal{H}ol(\mathbb{D})$, i.e. such that $R_{n}\circ b_{\lambda}-G\circ b_{\lambda}\in z^{n}\mathcal{H}ol(\mathbb{D})$, as $b_\lambda$ is self-inverse.

Due to invariance of the norm in $H^{\infty}$ with respect to the
composition by $b_{\lambda}$ we have 
\begin{align*}
\Norm{g}{H^{\infty}/b_{\lambda}^{n}H^{\infty}} & =\Norm{\Psi_{n}}{H^{\infty}/z^{n}H^{\infty}}\\
 & =\frac{1}{\abs{\zeta+r}}\Norm{\Phi_{n}}{H^{\infty}/z^{n}H^{\infty}},
\end{align*}
where
\[
\Phi_{n}(z)=1+r\left(-\frac{1+r\zeta}{\zeta+r}\right)^{n-1}z^{n}-\frac{1-r^{2}}{r+\zeta}\sum_{k=1}^{n-1}\left(-\frac{1+r\zeta}{r+\zeta}\right)^{k-1}z^{k},
\]
is an analytic polynomial of degree $n$ with nonnegative coefficients
because the condition $-\frac{1}{r}<\zeta<-1$ implies $-\frac{1+r\zeta}{r+\zeta}>0$.

Denote by $F_{n}$ the $(n-1)-$th Fejer kernel, $F_{n}(z)=\frac{1}{2\pi}\sum_{|j|\le n-1}\Big(1-\frac{|j|}{n}\Big)z^{j}$,
and denote by $*$ the usual convolution operation in $L^{1}(\partial\mathbb{D})$.
Then, for any $h\in L^{\infty}\left(\partial\mathbb{D}\right)$, we
have $\Norm{h*F_{n}}{\infty}\le\Norm{h}{\infty}\Norm{F_{n}}{H^{1}}=\Norm{h}{\infty}$.
On the other hand, since $\widehat{h_{1}*h_{2}}(j)=\hat{h_{1}}(j)\hat{h_{2}}(j)$
and $\widehat{F}_{n}(j)=0$ for every $j\ge n$, we have 
\[
h*F_{n}=\Phi_{n}*F_{n}
\]
for any $h\in H^{\infty}$ such that $\hat{h}(k)=\widehat{\Phi}_{n}(k)$,
$k=0,1,,\dots,n-1$. Hence, for any such $h$, $\Norm{h}{\infty}\ge\Norm{\Phi_{n}*F_{n}}{\infty}$
and so
\[
\begin{aligned} \Norm{\Phi_{n}}{H^{\infty}/z^{n}H^{\infty}}=\inf\big\{\|h\|_{\infty}:\:h\in H^{\infty},\:\hat{h}(k) & =\widehat{\Phi}_{n}(k),\:0\le k\le n-1\big\}\\
 & \geq\Norm{\Phi_{n}*F_{n}}{\infty}\ge(\Phi_{n}*F_{n})(1).
\end{aligned}
\]
Note that the convolution with $F_{n}$ gives us the Ces\`aro mean of
the partial sums of the Fourier series. Denoting by $S_{j}$ the $j-$th
partial sum for $\Phi_{n}$ at 1 we have
\[
(\Phi_{n}*F_{n})(1)=\frac{1}{n}\sum_{j=0}^{n-1}S_{j}(1).
\]
For $j=0, \dots, n-1$ we have 
\begin{align*}
S_{j}(1) & =1-\frac{1-r^{2}}{r+\zeta}\sum_{k=1}^{j}\left(-\frac{1+r\zeta}{r+\zeta}\right)^{k-1}\\
 & =\frac{r+\zeta+(1-r)\left(-\frac{1+r\zeta}{r+\zeta}\right)^{j}}{\zeta+1},
\end{align*}
and a computation shows that
\begin{align*}
\sum_{j=0}^{n-1}S_{j}(1) & =\frac{\zeta+r}{(\zeta+1)^{2}(1+r)} \left(\left(\frac{-r\zeta-1}{r+\zeta}\right)^{n}(r-1)+(-1+(\zeta+1)n)r+1+(\zeta+1)n\right).
\end{align*}
Passing to the limit as $\zeta\rightarrow-1$ we obtain
\begin{align*}
\lim_{\zeta\rightarrow-1}\Norm{(\zeta-M_{B})^{-1}}{} & =\Norm{(-1-M_{B})^{-1}}{}\\
 & =\lim_{\zeta\rightarrow-1}\Norm{g}{H^{\infty}/b_{-r}^{n}H^{\infty}}\\
 & =\lim_{\zeta\rightarrow-1}\Norm{\Psi_{n}}{H^{\infty}/z^{n}H^{\infty}}\\
 & \geq\frac{1}{1-r}\cdot\lim_{\zeta\rightarrow-1}\frac{1}{n}\sum_{j=0}^{n-1}S_{j}(1)\\
 & =\frac{n(1+r)+1-r}{2(1-r)},
\end{align*}
which completes the proof.
\end{proof}

\subsection{Proof of Theorem \ref{thm:Ptak_type_result}}

Of particular importance to our discussion in the proof of Theorem~\ref{thm:Ptak_type_result} will be the analytic Toeplitz
matrix $X_{\beta}$ given entry-wise by 
\begin{equation}
\left(X_{\beta}\right)_{ij}=\Bigg\{\begin{array}{ll}
0 & \mbox{if }i<j\\
1 & \mbox{if }i=j\\
\beta & \mbox{if }i>j
\end{array}\label{eq:X_xi_beta-1-1}
\end{equation}
with $\beta\in[0,2]$. The spectral norm of $X_{\beta}$ is computed
in \cite[Proposition II.6]{SzZa}. 
\begin{proof}
We first prove (\ref{eq:H_3}). We consider $T\in\cC_{n}$ with spectral
radius less than $r.$ It follows from \cite[Theorem III.2]{Sze}
that for any $\zeta\in\mathbb{C}-\sigma(T)$ the resolvent of $T$
is bounded by 
\begin{equation}
\Norm{R(\zeta,\,T)}{}\leq\Norm{R(\zeta,\,M_{B})}{},\label{eq:genius}
\end{equation}
where $B$ is the Blaschke product associated with $m_{T}$. We recall~\cite[Theorem III.2]{Sze}
that given a finite Blaschke product $B(z)=\prod_{i=1}^{\abs{m}}\frac{\lambda_{i}-z}{1-\bar{\lambda}_{i}z}$,
of degree $\abs{m}\geq1$, with zeros $\lambda_{i}\in\mathbb{D}$,
the components of the resolvent of the model operator $M_{B}$ at
any point $\zeta\in\mathbb{C}-\{\lambda_{1},...,\lambda_{\abs{m}}\}$
with respect to the Takenaka-Malmquist-Walsh basis are given by 
\begin{equation}
\left(\big(\zeta-M_{B}\big)^{-1}\right)_{1\leq i,\,j\leq\abs{m}}=\left\{ \begin{array}{ll}
0 & \mbox{if }i<j\\
\frac{1}{\zeta-\lambda_{i}} & \mbox{if }i=j\\
\frac{(1-\Abs{\lambda_{i}}^{2})^{1/2}}{\zeta-\lambda_{i}}\frac{(1-\Abs{\lambda_{j}}^{2})^{1/2}}{\zeta-\lambda_{j}}{\prod_{k=j+1}^{i-1}\frac{1-\bar{\lambda}_{k}\zeta}{\zeta-\lambda_{k}}} & \mbox{if }i>j
\end{array}.\right.\label{eq:OS_formula}
\end{equation}
 Recall that for any $n\times n$ matrices $A=(a_{ij})$ and $A=(a'_{ij})$,
the condition $\vert a_{ij}\vert\leq a'_{ij}$ implies that $\Norm{A}{}\leq\Norm{A'}{}$.
Hence, supposing $\vert\zeta\vert\geq1$ we have 
\[
\frac{1}{\left|\zeta-\lambda_{i}\right|}\leq\frac{1}{\left|\zeta\right|-\left|\lambda_{i}\right|}\leq\frac{1}{1-r}
\]
 on the one hand and 
\[
\frac{(1-\Abs{\lambda_{i}}^{2})^{1/2}}{\left|\zeta-\lambda_{i}\right|}\frac{(1-\Abs{\lambda_{j}}^{2})^{1/2}}{\left|\zeta-\lambda_{j}\right|}\leq\frac{(1-\Abs{\lambda_{i}}^{2})^{1/2}}{1-\left|\lambda_{i}\right|}\frac{(1-\Abs{\lambda_{j}}^{2})^{1/2}}{1-\left|\lambda_{j}\right|}\leq\frac{1+r}{1-r},
\]
for $\mbox{if }i>j$ on the other hand and estimate 
\[
\Abs{\big(\zeta-M_{B}\big)_{i,\,j}^{-1}}\leq\left\{ \begin{array}{ll}
0 & \mbox{if }i<j\\
\frac{1}{1-r} & \mbox{if }i=j\\
\frac{1+r}{1-r} & \mbox{if }i>j
\end{array}=\frac{1}{1-r}X_{1+r},\right.
\]
where the analytic Toeplitz matrix $X_{\beta}$ is entry-wise defined
by in (\ref{eq:X_xi_beta-1-1}). Moreover choosing $\lambda_{1}=\dots=\lambda_{n}=r$,
$T^{*}=M_{B}$ and $\zeta=1$ we get 
\begin{align*}
\left(\big(1-T^{*}\big)^{-1}\right)_{1\leq i,\,j\leq n} & =\left\{ \begin{array}{ll}
0 & \mbox{if }i<j\\
\frac{1}{1-r} & \mbox{if }i=j\\
\frac{1+r}{1-r} & \mbox{if }i>j
\end{array}\right.\\
 & =\frac{1}{1-r}X_{1+r}.
\end{align*}
This completes the proof of (\ref{eq:H_3}) and proves in particular
that for any $\abs{\zeta}\geq1$ 
\[
\Norm{R(\zeta,\,T)}{}\leq \Norm{\big(1-T^{*}\big)^{-1}}{}=\frac{1}{1-r}\Norm{X_{1+r}}{}.
\]
\end{proof}

\subsection{Proof of Theorem 2}
\begin{proof}
The spectral norm of $X_{\beta}$ is
computed in \cite[Proposition II.6]{SzZa}. In particular for any
$\beta\in[0,2]$ 
\[
\Norm{X_{\beta}}{}=\frac{1}{2}\sqrt{(\beta-2)^{2}+\frac{\beta^{2}}{\cot^{2}(\theta^{*}/2)}},
\]
where $\theta^{*}$ is the unique solution of 
\[
\cot(n\theta)=\frac{\beta-2}{\beta}\cot(\theta/2)
\]
in $[\frac{2n-1}{2n}\pi,\,\pi)$ and it follows, see~\cite{Sze},
that 
\[
\Norm{X_{0}}{}=1,\qquad\Norm{X_{1}}{}=\frac{1}{2\sin(\frac{\pi}{4n+2})},\qquad\Norm{X_{2}}{}=\cot \left(\frac{\pi}{4n}\right).
\]
In our situation $\beta=1+r$ and 
\[
\Norm{X_{1+r}}{}=\frac{1}{2}\sqrt{(1-r)^{2}+\frac{(1+r)^{2}}{\cot^{2}(\theta^{*}/2)}},
\]
where $\theta^{*}$ is the unique solution of 
\[
\cot(n\theta)=-\frac{1-r}{1+r}\cot(\theta/2)
\]
in $[\frac{2n-1}{2n}\pi,\,\pi).$ Note that if $\beta=1$ then $\theta^{*}=\frac{2n\pi}{2n+1}$
while if $\beta=2$ then $\theta^{*}=\frac{2n-1}{2n}\pi.$ We consider
the case $\beta=1+r$ and put 
\[
f_{n}(\theta)=\cot(n\theta)+\frac{1-r}{1+r}\cot(\theta/2).
\]
We first observe that
\begin{align*}
f_{n}\left(\frac{2n-1}{2n}\pi\right) & =\cot\left(\frac{2n-1}{2}\pi\right)+\frac{1-r}{1+r}\cot\left(\frac{2n-1}{4n}\pi\right)\\
 & =\frac{1-r}{1+r}\cot\left(\frac{2n-1}{4n}\pi\right)>0.
\end{align*}
Moreover $\lim_{\theta\rightarrow\pi^{-}}\cot(n\theta)=-\infty$ and
$\cot(\pi/2)=0$ so that
\[
\lim_{\theta\rightarrow\pi^{-}}f_{n}\left(\theta\right)=-\infty.
\]
The function $\theta\mapsto\cot(n\theta)$ is strictly decreasing
on $[(n-1)\pi/n,\pi)\supset[\frac{2n-1}{2n}\pi,\,\pi)$, and this
is also the case of the function $\theta\mapsto\frac{1-r}{1+r}\cot(\theta/2)$
on $(0,\pi)\supset[\frac{2n-1}{2n}\pi,\,\pi)$. Therefore the function
$f_{n}$ is also strictly decreasing on the interval $[\frac{2n-1}{2n}\pi,\,\pi)$,
which shows the existence and the uniqueness of $\theta^{*}\in[\frac{2n-1}{2n}\pi,\,\pi)$
such that
\[
f_{n}(\theta^{*})=0.
\]
Now we show that for $n$ large enough we have 
\[
f_{n}\left(\frac{2n}{2n+1}\pi\right)<0.
\]
On one hand we have
\begin{align*}
\cot\left(n\frac{2n}{2n+1}\pi\right) & =\cot\left(n\pi-\frac{\pi}{2}+\frac{\pi}{2(2n+1)}\right)\\
 & =-\tan\left(n\pi+\frac{\pi}{2(2n+1)}\right)\\
 & =-\frac{\sin\left(n\pi+\frac{\pi}{2(2n+1)}\right)}{\cos\left(n\pi+\frac{\pi}{2(2n+1)}\right)}\\
 & =-\frac{(-1)^{n}\sin\left(\frac{\pi}{2(2n+1)}\right)}{(-1)^{n}\cos\left(\frac{\pi}{2(2n+1)}\right)}\\
 & =-\tan\left(\frac{\pi}{2(2n+1)}\right)=-\frac{\pi}{4n}+\cO\left(\frac{1}{n^{2}}\right)
\end{align*}
as $n\rightarrow\infty$. On the other hand 
\[
\cot\left(\frac{n}{2n+1}\pi\right)=\frac{\pi}{4n}+\cO\left(\frac{1}{n^{2}}\right).
\]
Therefore 
\begin{align*}
f_{n}\left(\frac{2n}{2n+1}\pi\right) & =-\frac{\pi}{4n}\left(1-\frac{1-r}{1+r}\right)+\cO\left(\frac{1}{n^{2}}\right)\\
 & =-\frac{\pi r}{2(1+r)}\frac{1}{n}\left(1+\cO\left(\frac{1}{n}\right)\right)<0
\end{align*}
for $n$ large enough. Thus, for $\beta=1+r$ we have 
\[
\frac{2n-1}{2n}\pi\leq\theta^{*}\leq\frac{2n\pi}{2n+1},
\]
and 
\[
\cot(\theta^{*}/2)\sim\frac{\pi-\theta^{*}}{2}\sim\frac{\pi}{4n},
\]
as $n$ tends to $+\infty.$ In particular
\begin{align*}
\Norm{X_{1+r}}{} & \sim\frac{1}{2}\frac{1+r}{\cot(\theta^{*}/2)},\\
 & \sim\frac{2n(1+r)}{\pi}
\end{align*}
as $n$ tends to $+\infty.$ 
\end{proof}

\end{document}